\documentclass[10pt]{amsart} %%% i.e. use 12pt type
\usepackage{graphicx,latexsym,bm,amsmath,amssymb,verbatim,multicol,lscape}
\theoremstyle{definition}

\theoremstyle{remark}

\numberwithin{equation}{section}
\begin{document}
\title{The 2-adic valuations of Stirling numbers of the second kind}%
\author{Shaofang Hong}
\address{Yangtze Center of Mathematics, Sichuan University, Chengdu 610064, P.R. China}
%\curraddr{}
\email{sfhong@scu.edu.cn; s-f.hong@tom.com; hongsf02@yahoo.com}
\author{Jianrong Zhao}
\address{School of Economic Mathematics, Southwestern University of
Finance and Economics, Chengdu 610074, P.R. China}%
\email{mathzjr@foxmail.com}
\author{Wei Zhao}
\address{Science and Technology on Communication Security Laboratory, 
Chengdu 610041, P.R. China}%
\email{zhaowei9801@163.com}
\thanks{The research of Hong was supported partially by the National Science
Foundation of China Grant \# 10971145, by the Ph.D. Programs
Foundation of Ministry of Education of China Grant \#20100181110073 and by
Program for New Century Excellent Talents in University Grant \#
NCET-06-0785}

%\subjclass{Primary 11T22,11R18}%
\date{\today}%
%\dedicatory{}%
%\commby{}%
% ----------------------------------------------------------------
\begin{abstract}
In this paper, we investigate the $2$-adic valuations of the Stirling
numbers $S(n, k)$ of the second kind. We show that $v_2(S(4i,
5))=v_2(S(4i+3, 5))$ if and only if $i\not\equiv 7 \pmod {32}$. This
confirms a conjecture of Amdeberhan, Manna and Moll raised in 2008.
We show also that $v_2(S(2^n+1, k+1))= s_2(n)-1$ for any positive integer
$n$, where $s_2(n)$ is the sum of binary digits of $n$. It proves
another conjecture of Amdeberhan, Manna and Moll.
\end{abstract}
\maketitle
%----------------------------------------------------------------
\section{\bf Introduction}

Divisibility properties of integer sequences have long been objects
of interest in number theory. $p$-Adic valuation is the modern
language of divisibility. Given a prime $p$ and a positive integer
$m$, there exist unique integers $a$ and $n$, with $a$ not divisible
by $p$ and $n\ge 0$, such that $m=ap^n$. The number $n$ is called
the {\it $p$-adic valuation} of $m$, denoted by $n=v_p(m)$. Stirling
numbers are common topics in number theory and combinatorics. Let
$\mathbb{N}$ denote the set of natural numbers. The {\it Stirling
numbers of the first kind}, denoted by $s(n,k)$ (with a lower-case
"$s$"), count the number of permutations of $n$ elements with $k$
disjoint cycles. The {\it Stirling numbers of the second kind} $S(n,
k)$ (with a capital "$S$") is defined for $n\in \mathbb{ N}$ and
positive integer $k\le n$ as the number of ways to partition a set
of $n$ elements into exactly $k$ nonempty subsets. One can
characterize the Stirling numbers of the first kind and the Stirling
numbers of the second kind by
$$
(x)_n=\sum_{k=0}^ns(n,k)x^k \ {\rm and} \
x^n=\sum_{k=0}^nS(n,k)(x)_k
$$
respectively, where $(x)_n$ is the falling factorial
$(x)_n=x(x-1)(x-2)...(x-n+1)$. The Stirling numbers of the first and
second kind can be considered to be inverses of one another:
$$
\sum^{{\rm max}(j,k)}_{l=0}s(l,j)\cdot S(k,l)=\delta _{jk} \ {\rm
and} \ \sum^{{\rm max}(j,k)}_{l=0}S(l,j)\cdot s(k,l)=\delta _{jk},
$$
where $\delta _{jk}$ is the Kronecker delta. See \cite{[AMM1]}-\cite{[W]} for some
results on this topic.

Amdeberhan, Manna and Moll \cite{[AMM1]} studied the 2-adic
valuations of Stirling numbers of second kind. Actually, they
computed the 2-adic valuation $v_2(S(n, k))$ for $k\le 4$.
Furthermore, they found that the 2-adic valuation of the Stirling
numbers of second kind of order 5 is the first nontrivial case. In
fact, they showed that $v_2(S(4n+1, 5))=v_2(S(4n+2, 5))=0$ for all
natural numbers $n$. It was also observed in \cite{[AMM1]} that $v_2(S(4n,
5))=v_2(S(4n+3, 5))$ for most indices. Consequently, Amdeberhan,
Manna and Moll proposed a conjecture describing those indices $n$
such that $v_2(S(4n, 5))\ne v_2(S(4n+3, 5))$. That is,\\
\\
{\bf Conjecture 1.1.} \cite{[AMM1]} {\it $v_2(S(4n, 5))\ne v_2(S(4n+3, 5))$ if
and only if $n\in\{32j+7: j\in\mathbb{N}\}$. }\\

On the other hand, Lengyel \cite{[L]} conjectured, proved by
Wannemacker \cite{[W]}, a special case of the 2-adic valuation of
$S(n,k)$:
\begin{align*}
v_2(S(2^n,k))=s_2(k)-1,
\end{align*}
independently of $n$, where $s_2(k)$ means the base $2$ digital sum of $k$.
In \cite{[AMM1]}, Amdeberhan, Manna and Moll also posed the following conjecture.\\
\\
{\bf Conjecture 1.2}. \cite{[AMM1]} {\it  For all $k$ and $1\le k\le 2^n$, we
have}
$$
v_2(S(2^n+1,k+1))=s_2(k)-1.
$$

In this paper, we concern on the $2$-adic valuations of the Stirling
numbers of the second kind. We provide detailed analysis to the
$k$-level of the Stirling numbers of the second kind of order 5
where $5\le k\le 8$ (Note that this analysis was given in \cite{[AMM1]} when
$k=4$). Using this we then show that $v_2(S(4i, 5))\ne v_2(S(4i+3,
5))$ if and only if $i\in \{32j+7: j\ge 1\}$. This proves that
Conjecture 1.1 is true. We show also that $v_2(S(2^n+1,
k+1))=s_2(k)-1$ for any positive integer $n$, where $s_2(n)$ means
the sum of binary digits of $n$. This confirms Conjecture 1.2.

Throughout the paper, we will use several elementary properties of
$S(n,k)$, listed below:
\begin{align}\label{1}
S(n,k)=\frac{1}{k!}\sum_{j=1}^{k}(-1)^{k-j}{k\choose j}j^n.
\end{align}
% Relation Pochhammer
%\begin{align}
%x^n=\sum_{k=0}^nS(n,k)(x)_k,
%\end{align}
The generating function
\begin{align}
\frac{1}{(1-x)(1-2x)...(1-kx)}=\sum_{n=1}^{\infty}S(n,k)x^n.
\end{align}
The recurrence
\begin{align}\label{7}
S(n,k)=S(n-1,k-1)+kS(n-1,k).
\end{align}

\section{\bf Proof of Conjecture 1.1}

In the present section, we show the truth of Conjecture 1.1.
We begin with the following lemma.\\
\\
{\bf Lemma 2.1.}  {\it Let $N\ge 2$ be an
 integer and $r,i$ be odd numbers. For any  $m\in \mathbb{Z}^+$,
\begin{align}\label{0}
v_2((2^Nr\pm1)^{2^mi}-1)=m+N.
\end{align}
In particular, we have
$v_2(3^{2^m}-1)=v_2(5^{2^m}-1)=m+2$.}

\begin{proof}First we have
\begin{align}
(2^Nr\pm1)^{2^mi}-1=\pm2^{m+N}\cdot r\cdot
i+\sum_{j=2}^{2^mi}(\pm1)^j\cdot(2^Nr)^j\cdot{2^mi\choose
j}.\label{15}
\end{align}

On the other hand, for any integer $2\le j\le 2^mi$, since $N\ge2$,
we have
\begin{align}\label{16}
v_2((2^Nr)^j\cdot{2^mi\choose j})\ge Nj+m-v_2(j!)\ge Nj+m-j+1\ge
m+N+1.
\end{align}
Then it follows from (\ref{15}) and (\ref{16}) that (\ref{0}) is
true.
\end{proof}

Let $k\in \mathbb{N}$ be fixed and $m\in \mathbb{N}$. Then for $0\le
j <2^m$, we define
\begin{align}
 C_{m,j}:=\{2^mi+j: 2^mi+j\ge k \ {\rm and} \ i\in \mathbb{N}\}
\end{align}
and
\begin{align}
 v_2(C_{m,j}):=\{v_2(S(2^mi+j,k): 2^mi+j\ge k \ {\rm and} \ i\in \mathbb{N}\}.
\end{align}
Evidently, the classes $C_{m,j}$ form a partition of $\mathbb{N}$
into classes modulo $2^m$ and each class $C_{m,j}$ splits into
exactly two classes modulo $2^{m+1}$: $C_{m+1,j},C_{m+1,j+2^{m}}$.
The class $C_{m,j}$ is called {\it constant} if $ v_2(C_{m,j})$
consists of a singe value. This single value is called constant of
the class $ C_{m,j}$. Otherwise the class $C_{m,j}$ is called {\it
non-constant}. We can now define inductively $m$-{\it level}. The
{\it 1-level} consists of two classes:
$$C_{1,0}=\{2i: i\in \mathbb{N}\} \ {\rm and} \ C_{1,1}=\{2i+1: i\in \mathbb{N}\}.$$
Assume that for $m\ge 2$, the {\it $(m-1)$-level} has been defined
and that it consists of the $e$ classes
$$C_{m-1, i_1}, ..., C_{m-1, i_e}, \ {\rm where} \ 1\le i_1<...<i_e.$$
Each class $C_{m-1, i_j}$ splits into two classes modulo $2^m$, i.e. $C_{m,i_j}$ and
$C_{m,i_j+2^m}$. Then the {\it $m$-level} is formed by the nonconstant
classes modulo $2^m$.\\
\\
{\bf Lemma 2.2.} \cite{[AMM1]} {\it Each of the following is true:}

(i). $v_2(C_{4,4})=v_2(C_{4,7})=2$.

(ii). {\it $v_2(C_{4,12})>2$ and $v_2(C_{4,15})>2$.}

(iii). {\it The 4-level of $S(n,5)$ is
$\{C_{4,12},C_{4,15}\}$.}\\
\\
{\bf Lemma 2.3.}  {\it Each of the following is true:}

(i). $v_2(C_{5,12})=v_2(C_{5,15})=3.$

(ii). {\it $v_2(C_{5,28})>3$ and $v_2(C_{5,31})>3$.}

(iii). {\it The 5-level of $S(n,5)$ is $\{C_{5,28}, C_{5,31}\}$.}

\begin{proof}
By Lemma 2.2, we know that the 4-level of $S(n,5)$ is $\{C_{4,12},
C_{4,15}\}$ and  $v_2(C_{4,12})>2$ and $v_2(C_{4,15})>2$. It is
clear that $C_{4,12}=C_{5,12}\sqcup C_{5,28}$  and
$C_{4,15}=C_{5,15}\sqcup C_{5,31}$. Then we only need to consider
the four classes: $C_{5,12},C_{5,28},C_{5,15},$ $C_{5,31}$. By
(1.1), we know that the Stirling number $S(n,5)$ is given by
\begin{align}\label{2}
S(n,5)=\frac{1}{24}(5^{n-1}-4^n+2\cdot 3^n-2^{n+1}+1).
\end{align}
First we prove that (i) holds. By Lemma 2.1 we have $5^{2^5}\equiv
3^{2^5}\equiv 1\pmod {2^7}$. It then follows from (1.1) that for any
nonnegative integers $r$ and $t$, we have

\begin{align}
24S(2^5t+12,5)&=5^{2^5t+11}-4^{2^5t+12}+2\cdot
3^{2^5t+12}-2^{2^5t+13}+1\nonumber\\
\nonumber&\equiv 5^{11}+2\cdot3^{12}+1 \\
&\equiv 2^6 \pmod {2^7}\label{3}
\end{align}
and
\begin{align}
24S(2^5r+15,5)&=5^{2^5r+14}-4^{2^5r+15}+2\cdot
3^{2^5r+15}-2^{2^5r+16}+1\nonumber\\
\nonumber&\equiv 5^{14}+2\cdot3^{15}+1 \\
&\equiv 2^6 \pmod {2^7}.\label{4}
\end{align}
Thus we get
$$
S(2^5t+12,5)\equiv 8\pmod {2^4}$$ and
$$S(2^5r+15,5)\equiv 8\pmod {2^4}.
$$
These imply that $v_2(C_{5,12})=3$ and $v_2(C_{5,15})=3$. Part (i)
is proved.

Now we consider $C_{5,28},C_{5,31}$. By (\ref{2}) we have
\begin{align}
&24S(2^5t_1+28,5)-24S(2^5t_2+12,5)\nonumber\\
\nonumber&\equiv
5^{27}+2\cdot3^{28}-(5^{11}+2\cdot3^{12})\\
\nonumber&\equiv 5^{11}(5^{2^4}-1)+2\cdot 3^{12}(3^{2^4}-1)\\
\nonumber&\equiv 5^{11}2^6\alpha_1+2\cdot3^{12}2^6\alpha_2\\
&\equiv2^6 \pmod {2^7},\label{5}
\end{align}
where $t_1,t_2\in\mathbb{ N}$ and $\alpha_1$, $\alpha_2\in
1+2\mathbb{N}$. By (\ref{3}) and (\ref{5}) we get
$$24S(2^5t_1+28,5)\equiv0\pmod {2^7}.$$  Clearly $v_2(S(2^5t_1+28,5))>3$,
i.e., $v_2(C_{5,28})>3$. Similarly we compute that
\begin{align}
&24S(2^5r_1+31,5)-24S(2^5r_2+15,5)\nonumber\\
\nonumber&\equiv
5^{30}+2\cdot3^{31}-(5^{14}+2\cdot3^{15})\\
\nonumber&\equiv 5^{14}(5^{2^4}-1)+2\cdot 3^{15}(3^{2^4}-1)\\
\nonumber&\equiv 5^{14}2^6\alpha_3+2\cdot3^{15}2^6\alpha_4\\
&\equiv2^6 \pmod {2^7},\label{6}
\end{align}
where $r_1,r_2\in \mathbb{N}$ and $\alpha_3$, $\alpha_4\in
1+2\mathbb{N}$. Then by (\ref{4}) and (\ref{6}) we obtain that
$$24S(2^5r_2+31,5)\equiv0\pmod {2^7}.$$  Hence $v_2(C_{5,31})>3$.
Part (ii) is proved.

It is easy to check that $v_2(S(28,5))=6$, $v_2(S(60,5))=4$,
$v_2(S(31,5))=7$ and $v_2(S(61,5))=4$. Since $28,60\in C_{5,28}$ and
$31,63\in C_{5,31}$, $C_{5,28}$ and $C_{5,31}$ are non-constant
classes. The proof of Lemma 2.3 is complete.
\end{proof}

\noindent{\bf Lemma 2.4.} {\it Each of the following is true:}

(i). $v_2(C_{6,60})=v_2(C_{6,63})=4$.

(ii). {\it $v_2(C_{6,28})>4$ and $v_2(C_{6,31})>4$.}

(iii). {\it The 6-level of $S(n,5)$ is $\{C_{6,28}, C_{6,31}\}$.}

\begin{proof}
By Lemma 2.3, we have 5-level of $S(n,5)$ is $\{C_{5,28},
C_{5,31}\}$ and  $v_2(C_{5,28})>3$ and $v_2(C_{5,31})>3$. It is easy
to see that $C_{5,28}=C_{6,28}\sqcup C_{6,30}$  and
$C_{5,31}=C_{6,31}\sqcup C_{6,63}$. So it remains to deal with
$C_{6,28},C_{6,60},C_{6,31},$ $C_{6,63}$. First we consider
$C_{6, 28}$ and $C_{6, 60}$.  By Lemma 2.1 we have $5^{2^6}\equiv
3^{2^6}\equiv 1\pmod {2^8}$. Then
\begin{align*}
24S(2^6t+60,5)&=5^{2^6t+59}-4^{2^6t+60}+2\cdot
3^{2^6t+60}-2^{2^6t+61}+1\\
&\equiv 5^{59}+2\cdot3^{60}+1 \pmod {2^8},
\end{align*}
where $t\in \mathbb{N}$. Hence we conclude that
\begin{align*}
3^45^524S(2^6t+60,5)&\equiv 5^{2^6}3^4+2\cdot3^{2^6}5^5+3^45^5\\
&\equiv3^4+2\cdot5^5+3^45^5\\
&\equiv 2^7\pmod {2^8}.
\end{align*}
So $v_2(S(2^6t+60,5))=4$.  This implies that $v_2(C_{6,60})=4$.
Similarly we have
\begin{align*}
24S(2^6t+28,5)&=5^{2^6t+27}-4^{2^6t+28}+2\cdot
3^{2^6t+28}-2^{2^6t+29}+1\\
&\equiv 5^{27}+2\cdot3^{28}+1\\
&\equiv0 \pmod {2^8}.
\end{align*}
Thus $v_2(S(2^6t+28,5))>4$.  That is  $v_2(C_{6,28})>4$. We can
compute that $v_2(S(28,5))=6$ and $v_2(S(92,5))=5$. So $C_{6,28}$ is
non-constant.  In the same way, we obtain that $v_2(C_{6,63})=4$,
$v_2(C_{6,31})>4$ and $C_{6,31}$ is non-constant. This completes the
proof of Lemma 2.4.
\end{proof}

\noindent{\bf Lemma 2.5.}   {\it Each of the following is true:}

(i). $v_2(C_{7,92})=v_2(C_{7,95})=5$.

(ii). {\it $v_2(C_{7,28})>5$ and $v_2(C_{7,31})>5$.}

(iii). {\it The 7-level of $S(n,5)$ is $\{C_{7,28}, C_{7,31}\}$. }

\begin{proof}
By Lemma 2.4, we get  6-level of $S(n,5)$ is $\{C_{6,28},
C_{6,31}\}$ and  $v_2(C_{6,28})>4$ and $v_2(C_{6,31})>4$. Clearly
$C_{6,28}=C_{7,28}\sqcup C_{7,92}$  and $C_{6,31}=C_{7,31}\sqcup
C_{6,95}$.  First we consider $C_{7,28}$ and $C_{7,92}$. Since
$5^{2^7}\equiv 3^{2^7}\equiv 1\pmod {2^9}$ by Lemma 2.1,  then we
have
\begin{align*}
24S(2^7t+92,5)&=5^{2^7t+91}-4^{2^7t+91}+2\cdot
3^{2^7t+92}-2^{2^7t+93}+1\nonumber\\
\nonumber&\equiv 5^{91}+2\cdot3^{92}+1\\
&\equiv 2^8\pmod {2^9},
\end{align*}
where $t\in \mathbb{N}$. Hence $v_2(C_{7,92})=5$. Similarly
\begin{align*}
24S(2^7t+28,5)&=5^{2^7t+27}-4^{2^7t+28}+2\cdot
3^{2^7t+28}-2^{2^7t+29}+1\\
&\equiv 5^{27}+2\cdot3^{28}+1\\
&\equiv0 \pmod {2^9}.
\end{align*}
Thus  $v_2(C_{7,28})>5$. Since $v_2(S(28,5))=6$ and
$v_2(S(156,5))=11$,  $C_{7,28}$ is non-constant. In the same way, we
obtain that $v_2(C_{7,95})=5$, $v_2(C_{7,31})>5$ and $C_{7,31}$ is a
non-constant class. This completes the proof of Lemma 2.5.
\end{proof}

\noindent{\bf Lemma 2.6.} {\it Each of the following is true:}

(i). $v_2(C_{8,28})=v_2(C_{8,159})=6$.

(ii). {\it $v_2(C_{8,156})>6$ and $v_2(C_{8,31})>6$.}

(iii). {\it The 8-level of $S(n,5)$ is $\{C_{8,156}, C_{8,31}\}$.}

\begin{proof}
By Lemma 2.5, we  have 7-level of $S(n,5)$ is $\{C_{7,28},
C_{7,31}\}$, $v_2(C_{7,28})>5$ and $v_2(C_{7,31})>5$. Clearly
$C_{7,28}=C_{8,28}\sqcup C_{8,156}$  and $C_{7,31}=C_{8,31}\sqcup
C_{8,159}$.  Now we consider $C_{8,28}$ and $C_{8,156}$. Since
$5^{2^8}\equiv 3^{2^8}\equiv 1\pmod {2^{10}}$ by Lemma 2.1, we have
\begin{align*}
24S(2^8t+28,5)&=5^{2^8t+27}-4^{2^8t+28}+2\cdot
3^{2^8t+28}-2^{2^8t+29}+1\\
\nonumber&\equiv 5^{27}+2\cdot3^{28}+1\\
 \nonumber&\equiv 2^9\pmod {2^{10}},
\end{align*}
where $t\in \mathbb{N}$. Thus $v_2(C_{8,28})=6$. Similarly
\begin{align*}
24S(2^8t+156,5)&=5^{2^8t+155}-4^{2^8t+156}+2\cdot
3^{2^8t+156}-2^{2^8t+157}+1\\
\nonumber&\equiv 5^{155}+2\cdot3^{156}+1\\
\nonumber&\equiv 0 \pmod {2^{10}}.
\end{align*}
This implies that $v_2(C_{8,28})>6$. Since $v_2(S(156,5))=11$ and
$v_2(S(412,5))=7$,  $C_{8,156}$ is a non-constant class. In the same
way, we get $v_2(C_{8,159})=6$, $v_2(C_{8,31})>6$ and $C_{8,31}$ is
a non-constant class. This completes the proof of Lemma 2.6.
\end{proof}

We can now show that Conjecture 1.1 is true.\\
\\
{\bf Theorem 2.7.} {\it We have that $v_2(S(4i, 5))\ne v_2(S(4i+3,
5))$ if and only if  $i\in\{32j+7: j\in\mathbb{N}\}$.}

\begin{proof}
First let  $i\in\{32j+7:j\ge0\}$. Then  $4i=2^7j+28\in C_{7,28}$ and
$4i+3=2^7j+31\in C_{7,31}$ for some nonnegative integer $j$. Clearly
$C_{7,28}=C_{8,28}\sqcup C_{8,156}$ and $C_{7,31}=C_{8,31}\sqcup
C_{8,159}$. But by Lemma 2.6 (i) and (ii), we have
$v_2(C_{8,156})>6$, $v_2(C_{8,31})>6$ and
$v_2(C_{8,28})=v_2(C_{8,159})=6$. Then we obtain that
$$
v_2(S({2^8t+28},5))\neq v_2(S({2^8t+31},5))
$$
and
$$
v_2(S({2^8t+156},5))\neq v_2(S({2^8t+159},5))
$$
for any integer $t\ge 0$. This implies that
$$
v_2(S(2^7j+28,5))\neq v_2(S(2^7j+31,5)).
$$

Now we give a partition of $\{4i:i\in \mathbb{N}\}=C_{2,0}$:
$$C_{3,0},
C_{4,4},C_{5,12},C_{6,60},C_{7,92},C_{7,28}.
$$
We also get a partition of $\{4i+3:i\in \mathbb{N}\}=C_{2,3}$:
$$C_{3,3}, C_{4,7},C_{5,15},C_{6,63},C_{7,95},C_{7,31}.$$ By Lemmas
2.2-2.6, we have $v_2(C_{3,0})=v_2(C_{3,3})=1$,
$v_2(C_{4,4})=v_2(C_{4,7})=2$, $v_2(C_{5,12})=v_2(C_{5,15})=3$
,$v_2(C_{6,60})=v_2(C_{6,63})=4 $ and
$v_2(C_{7,92})=v_2(C_{7,95})=5$. Thus $v_2(S(4i, 5))=v_2(S(4i+3,
5))$ except that $i\in\{32j+7:j\ge0\}$. This completes the
proof of Theorem 2.7.
\end{proof}

By Theorem 2.7, we know immediately that $v_2(S(4i, 5))=v_2(S(4i+3,
5))$ if and only if $i\not\equiv 7 \pmod {32}$.

\section{\bf The 2-adic valuations of $S(2^n+1, k+1)$ and $S(2^n+2, k+2)$}

In this section, we investigate the 2-adic valuations of $S(2^n+1,
k+1)$ and $S(2^n+2, k+2)$. In 2005, Wannemacker computed
the 2-adic valuation of $S(2^n,k)$ as follows.\\
\\
{\bf Lemma 3.1.} \cite{[W]} {\it For all $k\in\mathbb{N}$ such that $1\le k\le 2^n$,
we have
\begin{align}
v_2(S(2^n,k))=s_2(k)-1.
\end{align} }

First we calculate the 2-adic valuation of $S(2^n+1, k+1)$. For this
purpose, we define $u(k)$ as follows: If $k=\sum_{i=0}^\infty
a_i2^i$, where $a_i\in\{0,1\}$, then $u(k):=0$ if $a_0=0$, and
$u(k):=t+1$ if $a_0=a_1=...=a_{t}=1$ and $a_{t+1}=0$. That is,
$u(k)$ is the first index $i$ such that $a_i=0$.\\
\\
{\bf Theorem 3.2.} {\it  For all $k\in\mathbb{N}$ such that $1\le k\le 2^n$, we have}
$$
v_2(S(2^n+1,k+1))=s_2(k)-1.
$$

\begin{proof}
Let $k=\sum_{i=0}^na_i2^i$, where $a_i\in\{0,1\}$. Then
$s_2(k)=\sum_{i=0}^{n}a_i$. We claim that
\begin{align}\label{8}
s_2(k+1)=s_2(k)+1-u(k).
\end{align}
If $u(k)=0$, then (3.2) is clearly true.
If $u(k)\ge 1$, then $a_0=1$. We may write
$k=1+2+...+2^{u(k)-1}+a_{u(k)+1}2^{u(k)+1}+...+a_n2^n$. Then
\begin{align}\label{9}
k+1=2^{u(k)}+a_{u(k)+1}2^{u(k)+1}+...+a_n2^n.
\end{align}
Thus $s_2(k)=u(k)+\sum_{i=u(k)+1}^n a_i$ and
$s_2(k+1)=1+\sum_{i=u(k)+1}^n a_i$. It infers that $s_2(k+1)-1=s_2(k)-u(k).$
The claim (3.2) is proved.

By (\ref{7}) and Lemma 3.1 we have
\begin{align}
\nonumber S(2^n+1,k+1)&=S(2^n,k)+(k+1)S(2^n,k+1)\\
&=2^{s_2(k)-1}\beta_1+(k+1)2^{s_2(k+1)-1}\beta_2,\label{10}
\end{align}
where $\beta_1$, $\beta_2\in 1+2\mathbb{N}$.

If $2|k$, then $u(k)=0$. Hence by (\ref{8}) and (\ref{10}) we get
\begin{align*}
S(2^n+1,k+1)=2^{s_2(k)-1}\beta_1+(k+1)2^{s_2(k)}\beta_2.
\end{align*}
It follows immediately that
$$v_2(S(2^n+1,k+1))=s_2(k)-1$$
as desired.

If $2\nmid k$, then $u(k)\ge 1$. By (\ref{8}) and (\ref{10}) we have
\begin{align*}
S(2^n+1,k+1)&=2^{s_2(k)-1}\beta_1+2^{u(k)}\beta\cdot
 2^{s_2(k)-u(k)}\beta_2\\
&=2^{s_2(k)-1}\beta_1+2^{s_2(k)}\beta\cdot\beta_2,
\end{align*}
where $k+1=2^{u(k)}\beta$ for some $\beta\in 1+2\mathbb{N}$ by (\ref{9}).
In this case we also have
$$v_2(S(2^n+1,k+1))=s_2(k)-1$$
as required. This completes the proof of Theorem 3.2.
\end{proof}

By Theorem 3.2, we know that Conjecture 1.2 is true. Finally, we
calculate the 2-adic valuation of $S(2^n+2,k+2)$.\\
\\
{\bf Theorem 3.3} {\it Let $k\in\mathbb{Z}$ and $1\le k\le 2^n$.

{\rm (i)}. If $u(k)=0$, then  $v_2(S(2^n+2,k+2))=s_2(k)-1.$

{\rm (ii)}. If $u(k)=1$, then $v_2(S(2^n+2,k+2))\ge s_2(k).$

{\rm (iii)}. If $u(k)\ge2$, then $v_2(S(2^n+2,k+2))=s_2(k)-u(k).$}

\begin{proof}
By (\ref{7}) we have
\begin{align}\label{11}
S(2^n+2,k+2)=S(2^n+1,k+1)+(k+2)S(2^n+1,k+2).
\end{align}
By Theorem 3.2  we get
\begin{align}\label{12}
v_2(S(2^n+1,k+1))=s_2(k)-1
\end{align}
and
\begin{align}\label{13}
v_2((k+2)S(2^n+1,k+2))=s_2(k+1)-1+v_2(k+2).
\end{align}

If $u(k)=0$, then $a_0=0$. By (\ref{8}) we get $s_2(k+1)-1=s_2(k).$
Thus by (\ref{11}), (\ref{12}) and (\ref{13}) we deduce that part
(i) is true.

If $u(k)\ge 1$, then $a_0=1$ and $v(k+2)=0$. Then by (\ref{8}) and
(\ref{13}) we have
\begin{align}\label{14}
v_2((k+2)S(2^n+1, k+2))=s_2(k)-u(k).
\end{align}
Hence by (\ref{11}), (\ref{12}) and (\ref{14}) we know that parts
(ii) and (iii) are true. The proof of Theorem 3.3 is complete.
\end{proof}

{\bf Acknowledgments}\\

The authors would like to thank the anonymous referee for his/her
helpful comments which improved its presentation.
\\

\end{document}